\documentclass[12pt,a4paper]{amsart}
\usepackage{geometry}
\usepackage{tabls}
\usepackage{url}
\usepackage[utf8]{inputenc}
\usepackage{amsfonts}
\usepackage{amssymb}
\usepackage{mathrsfs}
\usepackage{amsmath}
\usepackage{amsthm}
\usepackage{amscd}
\usepackage{fancyhdr}
\usepackage{enumerate}
\usepackage{paralist}
\usepackage{xypic}
\usepackage{graphicx}
\usepackage{cite}
\usepackage{lipsum}
\usepackage{footmisc}
\usepackage[all,cmtip]{xy}

\numberwithin{equation}{section}

\theoremstyle{plain}

\newtheorem{Thm}[equation]{Theorem}
\newtheorem{lem}[equation]{Lemma}
\newtheorem{prop}[equation]{Proposition}
\newtheorem{rem}[equation]{Remark}

\begin{document}

\title{Motivic double zeta values of odd weight}

\author{Jiangtao Li*}
\email{lijiangtao@amss.ac.cn}
\address{Jiangtao Li* \\Hua Loo-Keng Center for Mathematics Sciences,
          Academy of Mathematics and Systems Science,
         Chinese Academy of Sciences, 
         Beijing, China}

\author{Fei Liu}
\email{liufei54@pku.edu.cn}
\address{Fei Liu \\School of Mathematical Sciences,
        Peking University,
         Beijing, China}

\begin{abstract}

   For odd $N\geq 5$, we establish a short exact sequence about motivic double zeta values $\zeta^{\mathfrak{m}}(r,N-r)$ with $r\geq3$ odd, $N-r\geq2$. From this we classify all the relations among depth-graded motivic double zeta values $\zeta^{\mathfrak{m}}(r,N-r)$ with $r\geq3$ odd, $N-r\geq2$. As  a corollary, we confirm a conjecture of Zagier on the rank of a matrix which concerns relations among multiple zeta values of odd weight.
   \end{abstract}
\maketitle
\let\thefootnote\relax\footnotetext{
2010 $\mathnormal{Mathematics} \;\mathnormal{Subject}\;\mathnormal{Classification}$. Primary 11F32, Secondary 11F67.\\
$\mathnormal{Keywords:}$  Multiple zeta values, Period polynomial, mixed Tate Motives. }

\section{Introduction}

The multiple zeta values are defined by convergent series
\[\zeta(n_1,...,n_r)=\sum_{0<k_1<\cdots<k_r}\frac{1}{k_1^{n_1}\cdots k_r^{n_r}}, \, (n_1,\cdots,n_{r-1}>0,n_r>1).\]
We call $N=n_1+n_2+\cdots+n_r$ and $r$ the weight and depth of $\zeta(n_1,...,n_r)$ respectively.
Let $\mathcal{Z}_n$ be the set of $\mathbb{Q}$-linear combinations of multiple zeta values of weight $n$.
Denote $\mathcal{Z}_0=\mathbb{Q}$, then the graded vector space
\[\mathcal{Z}=\bigoplus\limits_{n\geq 0}^{+\infty}\mathcal{Z}_n\]
is a $\mathbb{Q}$-graded algebra. Beware that the notation about multiple zeta values here is different from that in \cite{ma}.

In order to study these numbers, Brown introduced the motivic multiple zeta algebra $\mathcal{H}$ in \cite{brown}. The elements of $\mathcal{H}$ are $\mathbb{Q}$-linear combinations of motivic multiple zeta values $\zeta^{\mathfrak{m}}(n_1,n_2,\cdots,n_r)$. We call $N=n_1+n_2+\cdots n_r$ and $r$  the weight and depth of $\zeta^{\mathfrak{m}}(n_1,n_2,\cdots,n_r)$ respectively. Denote by $\mathfrak{D}_r\mathcal{H}$ the  set of elements of $\mathcal{H}$ of depth $\leq r$ and $gr_r^{\mathfrak{D}}\mathcal{H}=\mathfrak{D}_r\mathcal{H}/\mathfrak{D}_{r-1}\mathcal{H}$. There is a natural graded $\mathbb{Q}$-algebra homomorphism
\[
\eta:\mathcal{H}\rightarrow\mathcal{Z},
\]
which satisfies $\eta(\zeta^{\mathfrak{m}}(n_1,n_2,\cdots,n_r))=\zeta(n_1,...,n_r)$.

Denote by $gr_r^{\mathfrak{D}}\mathcal{H}_N$ the weight $N$ part of $gr_r^{\mathfrak{D}}\mathcal{H}$, and denote by $\mathfrak{D}_N^{od,ev}$ the $\mathbb{Q}$-vector subspace of $gr_2^{\mathfrak{D}}\mathcal{H}$  generated by the natural images of $\zeta^{\mathfrak{m}}(r,N-r)$ for $3\leq r\leq N-2$ in $gr_2^{\mathfrak{D}}\mathcal{H}$. 

For even $k$,  denote by $\mathrm{W}^+_k$ and $\mathrm{W}^-_k$ the space of odd and even restricted period polynomials of weight $k$  respectively. Then we have
\begin{Thm}\label{exact}
For odd $N\geq 5$,  there is an exact sequence
\[
0\rightarrow \mathfrak{D}_N^{{od},{ev}}\xrightarrow{\partial} (gr_1^{\mathfrak{D}}\mathcal{H}^{od}\otimes_{\mathbb{Q}} gr_1^{\mathfrak{D}}\mathcal{H}^{ev})_N\xrightarrow{v}(W_{N-1}^+\oplus W_{N+1}^-)^{{}^{\vee}}\rightarrow 0.
\]
\end{Thm}

Here $(W_{N-1}^+\oplus W_{N+1}^-)^{{}^{\vee}}$ means the dual vector space of $W_{N-1}^+\oplus W_{N+1}^-$,
$gr_1^{\mathfrak{D}}\mathcal{H}^{od}$  is the subspace  of $gr_1^{\mathfrak{D}}\mathcal{H}$  generated by the images of  $$\zeta^{\mathfrak{m}}(3),\zeta^{\mathfrak{m}}(5),\cdots,\zeta^{\mathfrak{m}}({2n+1}),\cdots$$ in $gr_1^{\mathfrak{D}}\mathcal{H}$, $gr_1^{\mathfrak{D}}\mathcal{H}^{ev}$ is the subspace of $gr_1^{\mathfrak{D}}\mathcal{H}$  generated by the images of
$$\zeta^{\mathfrak{m}}(2),\zeta^{\mathfrak{m}}(4),\cdots,\zeta^{\mathfrak{m}}(2n),\cdots$$ in $gr_1^{\mathfrak{D}}\mathcal{H}$ and $(gr_1^{\mathfrak{D}}\mathcal{H}^{od}\otimes gr_1^{\mathfrak{D}}\mathcal{H}^{ev})_N$ means the weight $N$ part of $(gr_1^{\mathfrak{D}}\mathcal{H}^{od}\otimes_{\mathbb{Q}} gr_1^{\mathfrak{D}}\mathcal{H}^{ev})$.

The map $\partial$ comes from the motivic Galois action on $\mathcal{H}$.
The precise definitions of the maps $\partial$ and $ v$ will be given in Section $2$.  The most difficult part in the proof of Theorem \ref{exact} is the exactness at the middle $\mathrm{Im}\;\partial=\mathrm{Ker}\; v$.

 Denote by $\partial^{{}^{\vee}}$ and $ v^{{}^{\vee}}$ the dual map of $\partial$ and $v$ respectively. The strategy to prove $$\mathrm{Im}\;\partial=\mathrm{Ker}\; v$$ is the following: \\
 Step $1$: We prove that $\mathrm{Im}\; v^{{}^{\vee}}\subseteq \mathrm{Ker}\;\partial^{{}^{\vee}}$. This part is essentially the result of Zagier \cite{zag}.\\
 Step $2$: Notice that elements in $W_{N-1}^+$ are symmetric and elements in $W_{N+1}^-$ are anti-symmetric. We show that if $a\in \mathrm{Ker}\;\partial^{{}^{\vee}}$ satisfies the symmetry or anti-symmetry property, then $a\in \mathrm{Im}\;v^{{}^{\vee}}$.\\
 Step $3$: We show that if $a\in \mathrm{Ker}\;v^{{}^{\vee}}$, then its symmetrization $a^+$ and anti-symmetrization $a^-$ both belong to $\mathrm{Ker}\;v^{{}^{\vee}}$. Since
 \[
 a=\frac{1}{2}(a^++a^-),
 \]
 we have $a\in \mathrm{Im}\; v^{{}^{\vee}}$.

From the Step $2$ and Step $3$, we deduce that the injective map $(41)$ in \cite{zag} is actually an isomorphism.
Also Theorem \ref{exact} provides a motivic explanation for Zagier's result about the dimension of multiple double zeta values of odd weight.

For $N\geq 5$  odd, denote \[\mathfrak{R}_N=\{(c_r)_{3\leq r\leq N-2,\;\mathrm{odd}}\mid\sum_{r=1,\;\mathrm{odd}}^{N-2}c_r\zeta^{\mathfrak{m}}(r,N-r)\equiv 0 \;\mathrm{mod} \; \zeta^{\mathfrak{m}}(N),c_r\in\mathbb{Q}\}.\]
Since motivic multiple zeta values  satisfy the double shuffle relation \cite{souderes}, by the main results in Ma \cite{ma} there is an injective $\mathbb{Q}$-linear map
\[
\xi: W_{N-1}^+\oplus W_{N+1}^-\rightarrow\mathfrak{R}_N
\]
for odd $N$ defined as follows.
For \[p\in \mathrm{W}^+_{N-1},\;p(x+y,y)=\sum_{0\leq r \leq N}\binom{N-2}{r-1}b_{N-r,r}x^{N-r-1}y^{r-2},\]
\[
\xi(p)=(b_{N-r,r}-b_{r,N-r})_{3\leq r\leq N-2,\mathrm{odd}}.
\]
For \[p\in \mathrm{W}^-_{N+1},\;\frac{\partial}{\partial x}p(x+y,y)=\sum_{0\leq r \leq N}\binom{N-2}{r-1}c_{N-r,r}x^{N-r-1}y^{r-1},\]
\[
\xi(p)=(c_{N-r,r}-c_{r,N-r})_{3\leq r\leq N-2,\mathrm{odd}}.
\]

 The proof of the main theorem of  Ma \cite{ma} is based on the work of Gangl,  Kaneko and Zagier about the formal double zeta space. From Theorem \ref{exact}, we can deduce that
\begin{Thm}\label{isom}
The map $\xi: W_{N-1}^+\oplus W_{N+1}^-\rightarrow\mathfrak{R}_N$ is an isomorphism.

\end{Thm}
Thus for  $N$ odd, we can complete the classification of relations among depth-graded motivic double zeta values $\zeta^{\mathfrak{m}}(r,N-r)$ with  $r$ odd. Recall that for $N$ even, the classification of relations among formal double zeta values of weight $N$ is completed by the work of Baumard and  Schneps \cite{schneps}.

\section{Motivic Galois action}\label{galois}

First we will give a very short introduction to mixed Tate motives, the references are \cite{bf}, \cite{deligne}.
Denote by $\mathcal{MT}(\mathbb{Z})$ the category of  mixed Tate motives over $\mathbb{Z}$.
It is a Tannakian category with the natural fiber functor
$$\omega: \mathcal{MT}(\mathbb{Z})\rightarrow \mathrm{Vec}_{\mathbb{Q}}, M\mapsto \oplus \omega_{r}(M),$$
where
$$\omega_{r}(M)=\mathrm{Hom}_{\mathcal{MT}(\mathbb{Z})}(\mathbb{Q}(r), gr_{-2r}^{\omega}(M)).$$

Denote by $\pi_1(\mathcal{MT}(\mathbb{Z}))$ the Tannakian fundamental group of $\mathcal{MT}(\mathbb{Z})$ with respect to the fiber functor $\omega$.
Then
\[
  \pi_1(\mathcal{MT}(\mathbb{Z}))=\mathbb{G}_m \ltimes \mathnormal{U},
  \]
where $U$ is the pro-unipotent algebraic group with free Lie algebra generated by $\sigma_{2n+1}$ of degree $-(2n+1), n\geq1$.

Denote by ${}_0\Pi_1$ the motivic fundamental groupoid of $\mathbb{P}^1-\{0,1,\infty\}$ from $\overrightarrow{1}_0$ to $\overrightarrow{-1}_1$ (the tangent vector $\overrightarrow{1}$ at $0$ and the tangent vector $\overrightarrow{-1}$ at $1$ ). Its ring of regular functions over $\mathbb{Q}$ is 
\[
\mathcal{O}(_0\Pi_1)\cong\mathbb{Q}\langle e^0,e^1\rangle,
\]
where $\mathbb{Q}\langle e^0,e^1\rangle$ is equipped with the shuffle product.

For any word $ u_1 \cdots u_k$ in $e^0,e^1$,
if $\varepsilon,\eta\to 0$, it's easy to check that
$$\int\limits_{\varepsilon<t_{1}<\cdots <t_{k}<1-\eta}\omega_{u_{1}}(t_{1})\cdots \omega_{u_{k}}(t_{k})=P(log(\varepsilon), log(\eta))+O(\mathrm{sup}(\varepsilon|log(\varepsilon)|^{A}+\eta|log(\eta)|^{B})),$$
where $P$ is a polynomial over $\mathbb{R}$, $A$ and $B$ are  positive constants, $\omega_{e^0}(t)=\frac{dt}{t}$ and $\omega_{e^1}(t)=\frac{dt}{1-t}$. Define the $\mathbb{Q}$-linear map
\[
dch:\mathcal{O}(_0\Pi_1)\to \mathbb{R}
\]
by $dch(u_1\cdots u_k)=P(0,0)$. By a direct calculation one finds that the image of $dch$ coincides with the space $\mathcal{Z}$. From the theory of iterated integrals, the map $dch$ is an algebra homomorphism with respect to the shuffle product.

Since  $\mathcal{O}({}_0\Pi_1)$ is an inductive limit of mixed Tate motives, $\pi_1(\mathcal{MT}\left(\mathbb{Z})\right)$ has a natural action on $\mathcal{O}({}_0\Pi_1)$ (under the Tannakian correspondence $\mathcal{MT}(\mathbb{Z})\cong \mathrm{Rep}_{\mathbb{Q}}\pi_1(\mathcal{MT}\left(\mathbb{Z})\right)$). Thus $\mathcal{O}(U)$ has a natural coaction on $\mathcal{O}({}_0\Pi_1)$.

Let $I\subseteq \mathcal{O}({}_0\Pi_1)$ be the kernel of $dch$. Denote by $J^{\mathcal{MT}}\subseteq I$ the largest graded ideal  which is stable under the coaction of $\mathcal{O}(U)$. Define $\mathcal{H}=\mathcal{O}({}_0\Pi_1)/J^{\mathcal{MT}}$ and  define the motivic multiple zeta value $\zeta^{\mathfrak{m}}(n_1,\cdots,n_r)$ to be the image of the word
$$e^1(e^0)^{n_1-1}e^1\cdots e^1(e^0)^{n_r-1}$$ in $\mathcal{H}$.  

Define the depth filtration $\mathfrak{D}_r\mathcal{H}$  as the images of the $\mathbb{Q}$-linear combinations of
$$(e^0)^{i_0}e^1(e^0)^{i_1}e^1\cdots e^1(e^0)^{i_s},\;s\leq r$$  in $\mathcal{H}$.  Define
\[
gr_r^{\mathfrak{D}}\mathcal{H}=\mathfrak{D}_r\mathcal{H}/\mathfrak{D}_{r-1}\mathcal{H}.
\]

 Denote by ${}_x\Pi_y$ the motivic fundamental groupoid of $\mathbb{P}^1-\{0,1,\infty\}$ from $x$ to $y$, where $x,y\in \{\overrightarrow{1}_0,\overrightarrow{-1}_1\}$. We write $\overrightarrow{1}_0,\overrightarrow{-1}_1$ as $0,1$ respectively for short. Denote by $G$ the group of automorphisms of the groupoids ${}_x\Pi_y$ for $x,y\in \{0,1\}$ satisfying the following conditions:\\
 $(1)$ The action is compatible with the composition laws
 \[
 {}_x\Pi_y\times {}_y\Pi_z\to {}_x\Pi_z
 \]
 for all $x,y,z\in \{0,1\}$.\\
 $(2)$ The automorphism fixes the elements
 \[
 \mathrm{exp}(e_0)\in {}_0\Pi_0(\mathbb{Q}),\mathrm{exp}(e_1)\in {}_1\Pi_1(\mathbb{Q}).
 \]
 By Proposition $5.11$ in \cite{deligne}, ${}_x\Pi_y$ is a $G$-torsor.
 Since the action of $U$ on ${}_x\Pi_y$ for $x,y\in\{0,1\}$ respects $(1)$ and $(2)$, we have a natural morphism
 \[
 \phi: U\to G\cong{}_0\Pi_1.
 \]
 By the main results in \cite{brown}, $\phi$ is injective. Denote by $\mathfrak{g}$ the Lie algebra of ${U}$. We have an injective map
 \[
 i:\mathfrak{g}\to \mathrm{Lie}\;G\cong (\mathbb{L}(e_0,e_1),\{\,,\,\}).
 \]
 Here $(\mathbb{L}(e_0,e_1),\{\,,\,\})$ is the free Lie algebra generated by $e_0,e_1$ with the following Ihara Lie bracket
 \[
 \{f,g\}=[f,g]+D_f(g)-D_g(f)
 \]
 and $D_f$ is the derivation on $\mathbb{L}(e_0,e_1)$ which satisfies $D_f(e_0)=0$ and $D_f(e_1)=[e_1,f]$ for $f\in \mathbb{L}(e_0,e_1)$.

There is a natural decreasing depth filtration on $(\mathbb{L}(e_0,e_1),\{\,,\,\})$ defined by
\[
\mathfrak{D}^r\mathbb{L}(e_0,e_1)=\{\xi\in \mathbb{L}(e_0,e_1)|\,\mathrm{deg}_{e_1}\,\xi \geq r\}.
\]
By Theorem $6.8\,(i)$ in \cite{deligne}, we have
\[
i(\sigma_{2n+1})=\mathrm{ad}^{2n}(e_0)e_1+\mathrm{higher\;depth\;terms}.
\]

The action of $\mathfrak{g}$ on $\mathcal{H}$ factors through $(\mathbb{L}(e_0,e_1),\{\;,\,\})$.
Since the depth filtration on $\mathfrak{g}$ and the depth filtration on $\mathcal{H}$ are compatible, the element $\sigma_{2n+1}$ defines a canonical map (see Section $10$ in \cite{depth} for a similar discussion in the totally odd case)
\[
\partial_{2n+1}:gr_r^{\mathfrak{D}}\mathcal{H}\rightarrow gr_{r-1}^{\mathfrak{D}}\mathcal{H}.
\]

For odd $N\geq 5$, define a map
\[
\partial:\mathfrak{D}_N^{od,ev}\rightarrow (gr_1^{\mathfrak{D}}\mathcal{H}^{od}\otimes_{\mathbb{Q}} gr_1^{\mathfrak{D}}\mathcal{H}^{ev})_N,
\]
\[
\partial(\zeta^{\mathfrak{m}}(r,N-r))=\sum_{3\leq s\leq N-2,odd}\zeta^{\mathfrak{m}}(s)\otimes\partial_s(\zeta^{\mathfrak{m}}(r,N-r)).
\]
By Theorem 3.3 in \cite{brown}, we know

\begin{prop}\label{injec}
For odd $N\geq 5$, $\partial$ is injective.
\end{prop}

We now give an explicit formula for $\partial\left(\zeta^{\mathfrak{m}}(2m+1,2n) \right)$. This can be achieved by the motivic method developed by Brown in \cite{depth}. Alternatively, we use a motivic lift of Zagier's explicit formula for $\zeta(2m+1,2n)$ to make the proof short.

\begin{prop}\label{matrix}
Suppose that $N\geq 5$ is odd and $m,n\geq1$ satisfy $2m+2n+1=N$,  then we have
\[
\begin{split}
\partial(\zeta^{\mathfrak{m}}&(2m+1,2n))=\\
&\sum_{\substack{m_1+n_1=\frac{1}{2}N-\frac{1}{2}\\
  m_1,n_1\geq1}}\left[\delta\dbinom{m_1,n_1}{m,n}-\binom{2m_1}{2m}-\binom{2m_1}{2n-1}\right]\zeta^{\mathfrak{m}}(2m_1+1)\otimes\zeta^{\mathfrak{m}}(2n_1),
\end{split}
\]
where $\delta\binom{m_1,n_1}{m,n}=1$  if $(m_1,n_1)=(m,n)$, $\delta\binom{m_1,n_1}{m,n}=0$  if $(m_1,n_1)\neq(m,n)$ and $\binom{n}{k}=\frac{n(n-1)\cdots (n-(k-1))}{k(k-1)\cdots 1}$ for $k\geq 1$.
\end{prop}
\noindent{\bf Proof}:
Since motivic multiple zeta values satisfy the double shuffle relation, and the proof of Proposition $7$ in \cite{zag} only uses the double shuffle relation among multiple zeta values, it follows from Proposition $7$ in \cite{zag} that for odd $N\geq 5$, $2m+2n+1=N, m,n\geq1$, we have
\[
\begin{split}
\zeta^{\mathfrak{m}}&(2m+1,2n)=-\frac{1}{2}\left[-\binom{N-1}{2n-1}-\binom{N-1}{2m}+1\right]\zeta^{\mathfrak{m}}(N)                        \\
                    &+\sum_{\substack{m_1+n_1=\frac{1}{2}N-\frac{1}{2}\\
  m_1,n_1\geq1}}\left[\delta\dbinom{m_1,n_1}{m,n}-\binom{2m_1}{2m}-\binom{2m_1}{2n-1}\right]\zeta^{\mathfrak{m}}(2m_1+1)\zeta^{\mathfrak{m}}(2n_1).
\end{split}
\]

  One computes
\[
\begin{split}
&\;\;\partial(\zeta^{\mathfrak{m}}(2m_1+1)\zeta^{\mathfrak{m}}(2n_1))\\
&=\sum_{1\leq k\leq\frac{N}{2}-\frac{3}{2}}\zeta^{\mathfrak{m}}(2k+1)\otimes \partial_{2k+1}(\zeta^{\mathfrak{m}}(2m_1+1)\zeta^{\mathfrak{m}}(2n_1))\\
&=\sum_{1\leq k\leq\frac{N}{2}-\frac{3}{2}}\zeta^{\mathfrak{m}}(2k+1)\otimes
\left(\partial_{2k+1}(\zeta^{\mathfrak{m}}(2m_1+1))\zeta^{\mathfrak{m}}(2n_1)+\zeta^{\mathfrak{m}}(2m_1+1)\partial_{2k+1}(\zeta^{\mathfrak{m}}(2n_1))\right)\\
&=\zeta^{\mathfrak{m}}(2m_1+1)\otimes\zeta^{\mathfrak{m}}(2n_1),
\end{split}
\]
where for the second equation we have used the fact $\partial$ is a derivative and for the last equation we have used the identities
\[
\partial_{2k+1}(\zeta^{\mathfrak{m}}(2m_1+1))=\delta\binom{2m_1+1}{2k+1},\partial_{2k+1}(\zeta^{\mathfrak{m}}(2n_1))=0.
\]
Here $\delta\binom{m}{n}$ is the usual Kronecker delta, namely, $\delta\binom{m}{n}=1$ if $m=n$, and $\delta\binom{m}{n}=0$ if $m\neq n.$

So we have
\[
\begin{split}
\partial(\zeta^{\mathfrak{m}}&(2m+1,2n))=\\
&\sum_{\substack{m_1+n_1=\frac{1}{2}N-\frac{1}{2}\\
  m_1,n_1\geq1}}\left[\delta\dbinom{m_1,n_1}{m,n}-\binom{2m_1}{2m}-\binom{2m_1}{2n-1}\right]\zeta^{\mathfrak{m}}(2m_1+1)\otimes\zeta^{\mathfrak{m}}(2n_1).
\end{split}
\]

$\hfill\Box$\\

   For an even integer $h>2$, if $P\in\mathbb{Q}[X]$ satisfies
\[
\begin{split}
                 &P(X)+X^{h-2}P(\frac{-1}{X})=0,\\
P(X)+&X^{h-2}P(1-\frac{1}{X})+(X-1)^{h-2}P(\frac{-1}{X-1})=0,
\end{split}
\]
then $P$ is called a period polynomial of weight $h$ over $\mathbb{Q}$.
Denote by $W_h$ the space of period polynomials of weight $h$ over $\mathbb{Q}$.

Define
\[
W_h^+=\{P\in W_h\mid P(0)=0, P(X)=X^{h-2}P(\frac{1}{X})\},
\]
\[
W_h^-=\{P\in W_h\mid P(0)=0, P(X)+X^{h-2}P(\frac{1}{X})=0\}.
\]
By the Eichler-Shimura-Manin correspondence, we have $$W_h=W_h^+\oplus W_h^-\oplus\mathbb{Q}(X^{h-2}-1).$$
For odd $N\geq 5$, define
\[
j_1:W^+_{N-1}\rightarrow ((gr_1^{\mathfrak{D}}\mathcal{H}^{od}\otimes gr_1^{\mathfrak{D}}\mathcal{H}^{ev})_N)^{{}^{\vee}}
\]
 for $p=\sum_{\substack{r+s=N-3\\r,s\geq1,odd}}p_{r,s}X^{r}\in W_{N-1}^+$ by
\[j_1(p)(\zeta^{\mathfrak{m}}(2m+1)\otimes\zeta^{\mathfrak{m}}(2n))=p_{2m-1,2n-1}.\]
Define
\[
j_2:W^-_{N+1}\rightarrow ((gr_1^{\mathfrak{D}}\mathcal{H}^{od}\otimes gr_1^{\mathfrak{D}}\mathcal{H}^{ev})_N)^{{}^{\vee}}
\]
 for $$q\in W_{N+1}^-\;\mathrm{ and}\; X^{N-2}q\,'(\frac{1}{X})=\sum_{\substack{r+s=N-1\\r,s\geq 1,odd}}t_{r-1,s-1}X^{r-1}$$ by
\[
j_2(q)(\zeta^{\mathfrak{m}}(2m+1)\otimes\zeta^{\mathfrak{m}}(2n))=t_{2m,2n},
\]
where $q\,'(X)$ means the derivative of $q(X)$.

   It is clear that $j_1$ and $j_2$ are both injective. We define
\[
j:W^+_{N-1}\oplus W^-_{N+1}\rightarrow ((gr_1^{\mathfrak{D}}\mathcal{H}^{od}\otimes gr_1^{\mathfrak{D}}\mathcal{H}^{ev})_N)^{{}^{\vee}}
\]
to be the unique linear map which satisfies $j\mid_{W^+_{N-1}}=j_1$ and $j\mid_{W^-_{N+1}}=j_2$. By the symmetry property of $W^+_{N-1}$ and the anti-symmetry property of $W^-_{N+1}$ it is easy to show that $j$ is also injective.

Define $v:(gr_1^{\mathfrak{D}}\mathcal{H}^{od}\otimes gr_1^{\mathfrak{D}}\mathcal{H}^{ev})_N\rightarrow (W^+_{N-1}\oplus W^-_{N+1})^{{}^{\vee}}$ to be the dual of $j$. So we have

\begin{prop}\label{surj}
For odd $N\geq 5$, the map
\[
v:(gr_1^{\mathfrak{D}}\mathcal{H}^{od}\otimes gr_1^{\mathfrak{D}}\mathcal{H}^{ev})_N\rightarrow (W^+_{N-1}\oplus W^-_{N+1})^{{}^{\vee}}
\]
is surjective.
\end{prop}

In the rest of this section, we prove
\begin{prop}\label{contain}
For odd $N\geq 5$, we have $\mathrm{Im}\;\partial\subseteq \mathrm{Ker}\;v$.
\end{prop}
\noindent{\bf Proof}:
By  the definition of $v$, it suffices to prove that
\[
\mathrm{Im}\;j\subseteq \mathrm{Ker}\; \partial^{{}^{\vee}}.
\]
From the explicit formula of $\partial$, this is equivalent to \\
(1) If
\[
p(X)=\sum_{\substack{m_1+n_1=\frac{1}{2}N-\frac{1}{2}\\m_1,n_1\geq 1}}p_{2m_1-1,2n_1-1}X^{2m_1-1}\in W_{N-1}^+,
\]
then
\[
\sum_{\substack{m_1+n_1=\frac{1}{2}N-\frac{1}{2}\\m_1,n_1\geq 1}}p_{2m_1-1,2n_1-1}\left[-\delta\binom{m_1,n_1}{m,n}+\binom{2m_1}{2m}+\binom{2m_1}{2n-1}\right]=0
\]
for all $m,n\geq1$ with $m+n=\frac{1}{2}N-\frac{1}{2}$.\\
(2) If $q\in W_{N+1}^-$ and
\[
X^{N-2}q\,'(\frac{1}{X})=\sum_{\substack{m_1+n_1=\frac{1}{2}N-\frac{1}{2}\\m_1,n_1\geq 1}}t_{2m_1,2n_1}X^{2m_1},
\]
then
\[
\sum_{\substack{m_1+n_1=\frac{1}{2}N-\frac{1}{2}\\m_1,n_1\geq 1}}t_{2m_1,2n_1}\left[-\delta\binom{m_1,n_1}{m,n}+\binom{2m_1}{2m}+\binom{2m_1}{2n-1}\right]=0
\]
for all $m,n\geq1$ with $m+n=\frac{1}{2}N-\frac{1}{2}$.\\
The statements $(1)$ and $(2)$ follow from the main results of Zagier in Section $6$ \cite{zag}. $\hfill\Box$\\

\section{Exactness at the middle}\label{key}

 In this section, we fix  an odd number  $N \geq 5$. From  the definition of $v$, $\mathrm{Im}\;\partial=\mathrm{Ker}\;v$
is equivalent to $\mathrm{Im}\;j=\mathrm{Ker}\;\partial^{{}^{\vee}}$.
By the explicit formula of $\partial^{{\vee}}$, $\mathrm{Im}\;j=\mathrm{Ker}\;\partial^{{}^{\vee}}$ is equivalent to
\begin{prop}\label{poly}
If $C(X)=\sum\limits_{n=1}^{\frac{1}{2}(N-3)}c_{n} X^{2n}$ satisfies
\[
C(X)-C(1+X)-X^{N-2}C(1+\frac{1}{X})=c+\mathrm{odd\;polynomial\;in}\;X,
\]
then there exist $p\in W_{N-1}^+$ and $q\in W_{N+1}^-$ such that $C(X)=Xp(X)+X^{N-2}q\,'(\frac{1}{X})$.
\end{prop}

For convenience, we write $L_C(X)=C(X)-C(1+X)-X^{N-2}C(1+\frac{1}{X})$.
Notice that elements in $W_{N-1}^+$ are symmetric and elements in $W_{N+1}^-$ are anti-symmetric. We first prove that if $C(X)$ satisfies some kinds of symmetry and anti-symmetry properties, then $C(X)$ comes from an odd and even period polynomial, respectively.
\begin{lem}\label{odd}
Suppose that $C(X)=\sum\limits_{n=1}^{\frac{1}{2}(N-3)}c_{2n} X^{2n}$ satisfies
\[
C(X)-C(1+X)-X^{N-2}C(1+\frac{1}{X})=c+\mathrm{odd\;polynomial\;in}\;X
\]
and $p(X):=\frac{C(X)}{X}$ satisfies $p(X)=X^{N-3}p(\frac{1}{X})$, then we have $p\in W_{N-1}^+$.
\end{lem}
\noindent{\bf Proof}:
By assumption, we have
\[
(O)\;\;\;\;(X+1)[p(X)-p(1+X)-X^{N-3}p(1+\frac{1}{X})]=L(X)+p(X)=c+\mathrm{odd\;polynomial\;in}\;X.
\]
It is easily seen that $c=L_{C}(0)=-C(1)$.
Replacing $X$ by $\frac{1}{X}$, we get
\[
(\frac{1}{X}+1)[p(\frac{1}{X})-p(1+\frac{1}{X})-\frac{1}{X^{N-3}}p(1+X)]=L(\frac{1}{X})+p(\frac{1}{X})=c+\mathrm{odd\;polynomial\;in}\;\frac{1}{X}.
\]
Multiplying the above by $X^{N-2}$ and using the assumption $p(X)=X^{N-3}p(\frac{1}{X})$, we get
\[
(E)\;\;\;\;(X+1)[p(X)-p(1+X)-X^{N-3}p(1+\frac{1}{X})]=cX^{N-2}+\mathrm{even\;polynomial\;in}\;X.
\]
By comparing formulas $(O)$ and $(E)$ , we have
\[
\;p(X)-p(1+X)-X^{N-3}p(1+\frac{1}{X})=c\frac{(1+X^{N-2})}{1+X}.
\]
Letting $X\rightarrow-1$ in the above formula, we obtain \[-p(1)=-C(1)=c(N-2)=-C(1)(N-2).\]
So $c=-C(1)=0$ and $p\in W_{N-1}^+$.
$\hfill\Box$\\

\begin{lem}\label{ev}
Suppose that $C(X)=\sum\limits_{n=1}^{\frac{1}{2}(N-3)}c_{n} X^{2n}$ satisfies
\[
C(X)-C(1+X)-X^{N-2}C(1+\frac{1}{X})=c+\mathrm{odd\;polynomial\;in}\;X
\]
and $q(X)$ satisfies $C(X)=X^{N-2}q\,'(\frac{1}{X})$, $q(X)+X^{N-1}q(\frac{1}{X})=0$ and $q(0)=0$, then we have $q\in W_{N+1}^-$.
\end{lem}
\noindent{\bf Proof}:
By assumption, we have
\[
\begin{split}
(O\,')\;\;\;\;(X+1)^{N}&\frac{d}{dX}\left[q(\frac{1}{X+1})-q(\frac{X}{X+1})+(X+1)^{-N+1}q(X)\right]\\
                          &=L(X)+q\,'(X)= c+\mathrm{odd\;polynomial\;in}\;X.
\end{split}
\]
Replacing $X$ by $\frac{1}{X}$, we get
\[
\begin{split}
\;\;\;\;\;\;-(\frac{1}{X}+1)^{N}&X^2\frac{d}{dX}\left[q(\frac{X}{X+1})-q(\frac{1}{X+1})+(\frac{1}{X}+1)^{-N+1}q(\frac{1}{X})\right]\\
                          &=L(\frac{1}{X})+q\,'(\frac{1}{X})= c+\mathrm{odd\;polynomial\;in}\;\frac{1}{X}.
\end{split}
\]
Multiplying the above by $X^{N-2}$, we get
\[
\begin{split}
(E\,')\;\;\;\;(X+1)^{N}&\frac{d}{dX}\left[q(\frac{1}{X+1})-q(\frac{X}{X+1})+(X+1)^{-N+1}q(X)\right]\\
                          &= cX^{N-2}+\mathrm{even\;polynomial\;in}\;X.
\end{split}
\]
By comparing formulas $(O\;')$ and $(E\;')$, we find
\[
\frac{d}{dX}\left[q(\frac{1}{X+1})-q(\frac{X}{X+1})+(X+1)^{-N+1}q(X)\right]=c\frac{(X^{N-2}+1)}{(X+1)^N}.
\]
Taking integration on both sides, we get
\[
q(\frac{1}{X+1})-q(\frac{X}{X+1})+(X+1)^{-N+1}q(X)=\frac{c}{N-1}\frac{X^{N-1}-1}{(X+1)^{N-1}}+a.
\]
Replacing $X$ by $\frac{1}{X}$ and using the assumption $q(X)+X^{N-1}q(\frac{1}{X})=0$, we get
\[
-q(\frac{1}{X+1})+q(\frac{X}{X+1})-(X+1)^{-N+1}q(X)=-\frac{c}{N-1}\frac{X^{N-1}-1}{(X+1)^{N-1}}+a.
\]
So $a=0$ and the above formula reduces to
\[
-q(\frac{1}{X+1})+q(\frac{X}{X+1})-(X+1)^{-N+1}q(X)=-\frac{c}{N-1}\frac{X^{N-1}-1}{(X+1)^{N-1}}.
\]
Letting $X=0$ in the above formula, we obtain $c=-(N-1)q(1)=0$.

So $q(\frac{1}{X+1})-q(\frac{X}{X+1})+(X+1)^{-N+1}q(X)=0$, $q\in W_{N+1}^-$.  $\hfill\Box$\\

Lemma \ref{odd} and \ref{ev} show that if $C(X)$ satisfies some kind of symmetry or anti-symmetry property, then Proposition \ref{poly} is true. The final step is to show that  there is some kind of symmetrization process.
\begin{lem}\label{sym}
Suppose that $C(X)=\sum\limits_{n=1}^{\frac{1}{2}N-\frac{3}{2}}c_{n}X^{2n}$ satisfies
\[
C(X)-C(1+X)-X^{N-2}C(1+\frac{1}{X})=\mathrm{constant}+\mathrm{odd\;polynomial\;in\;}X.
\]
  Then the polynomial    $p(X)=C\;'(X)+X^{N-3}C\;'(\frac{1}{X})$  satisfies
\[
Xp(X)-(1+X)p(1+X)-X^{N-2}(1+\frac{1}{X})p(1+\frac{1}{X})=\mathrm{constant}+\mathrm{odd\;polynomial\;in\;}X.
\]
\end{lem}
\noindent{\bf Proof}:
Let $N=2K+1$, it suffices to prove that
\[
\sum_{j=1}^{K-1}(jc_j+(K-j)c_{K-j})\left[\binom{2j}{2i}+\binom{2j}{2K-2i-1}-\delta\binom{j}{i}\right]=0
\]
for $1\leq i\leq K-1$,
i.e.
\[
\begin{split}
\sum_{j=1}^{K-1}jc_j&\left[\binom{2j}{2i}+\binom{2j}{2K-2i-1}-\delta\binom{j}{i}\right]+\\
                    &\sum_{j=1}^{K-1}jc_j\left[\binom{2K-2j}{2i}+\binom{2K-2j}{2K-2i-1}-\delta\binom{K-j}{i}\right]=0
\end{split}
\]
for $1\leq i\leq K-1$, where $\delta\binom{m}{n}=1$ if $m=n$, and $\delta\binom{m}{n}=0$ if $m\neq n$.

First we fix some notations. For an indeterminate $T$ and an integer $n\geq1$, denote 
\[
\binom{T}{n}=\frac{T(T-1)\dots(T-n+1)}{n!}\in \mathbb{Q}[T]
\]
and for $n=0$ denote $\binom{T}{0}=1$.
For integers $n<0,m>0$, define $\binom{m}{n}=0$.

It is well-known that the polynomials
\[
\binom{T}{0},\binom{T}{1},\cdots,\binom{T}{n}
\]
form a basis of the vector space of polynomials of degree $\leq n$.
In fact
\[
f(0)\binom{T}{0}+\Delta f(0)\binom{T}{1}+\cdots+\Delta^n f(0)\binom{T}{n}=f(T),
\]
where $\Delta^lf(0)$ is the $l$-th difference of $f$ at $T=0$.

The differences are defined as follows:
\[\Delta^0f(T)=f(T),\]
\[\Delta^1f(T)=f(T+1)-f(T),\]
\[\Delta^2f(T)=\Delta^1f(T+1)-\Delta^1f(T)=f(T+2)-2f(T+1)+f(T),\]
\[\cdots\;\cdots\;\cdots\]

Let  $f(T)=\frac{T}{2}\left[\binom{T}{2i}+\binom{2K-T}{2i}\right]$,  we have
\[
\begin{split}
\Delta^1 f(T)&=f(T+1)-f(T)\\
           &=\frac{T+1}{2}\left[\binom{T+1}{2i}+\binom{2K-T-1}{2i}\right]-\frac{T}{2}\left[\binom{T}{2i}+\binom{2K-T}{2i}\right]\\
        &=\frac{T}{2}\binom{T}{2i-1}+\frac{1}{2}\binom{T+1}{2i}-\frac{T}{2}\binom{2K-T-1}{2i-1}+\frac{1}{2}\binom{2K-T-1}{2i},
\end{split}
\]
\[
\begin{split}
&\Delta^2 f(T)=\frac{T+1}{2}\binom{T+1}{2i-1}-\frac{T}{2}\binom{T}{2i-1}+\frac{1}{2}\left[\binom{T+2}{2i}-\binom{T+1}{2i}\right]\\
&-\left[\frac{T+1}{2}\binom{2K-T-2}{2i-1}-\frac{T}{2}\binom{2K-T-1}{2i-1}\right]+\frac{1}{2}\left[\binom{2K-T-2}{2i}-\binom{2K-T-1}{2i}\right]\\
&=\frac{T}{2}\binom{T}{2i-2}+\frac{2}{2}\binom{T+1}{2i-1}+\frac{T}{2}\binom{2K-T-2}{2i-2}-\frac{2}{2}\binom{2K-T-2}{2i-1}.
\end{split}
\]
For $0\leq p\leq 2i+1$, by induction we have
\[
\Delta^pf(T)=\frac{T}{2}\binom{T}{2i-p}+\frac{p}{2}\binom{T+1}{2i-p+1}+(-1)^p\frac{T}{2}\binom{2K-T-p}{2i-p}-(-1)^p\frac{p}{2}\binom{2K-T-p}{2i-p+1}.
\]

For $0\leq p\leq 2i+1$, denote  \[x_p=\Delta^pf(0)=\frac{p}{2}\left[\binom{1}{2i-p+1}-(-1)^p\binom{2K-p}{2i-p+1}  \right].\]
For $p> 2i+1$, the above formula gives $x_p=0$. We have
\[
(a)\;\;\;\;\;\;\;\sum_{p=1}^{2i+1}x_p\binom{T}{p}=\frac{T}{2}\left[\binom{T}{2i}+\binom{2K-T}{2i}\right].
\]
For $l\geq1$, by taking the $l$-th difference of the above formula on both sides, we get
\[
\begin{split}
&(b)\;\;\;\;\;\sum_{p=0}^{2i+1}x_{p+l}\binom{T}{p}\\
&=\frac{T}{2}\binom{T}{2i-l}+\frac{l}{2}\binom{T+1}{2i-l+1}+(-1)^l\frac{T}{2}\binom{2K-T-l}{2i-l}-(-1)^l\frac{l}{2}\binom{2K-T-l}{2i-l+1}.
\end{split}
\]

For $1\leq j\leq K-1$, by letting $T=2j$ in formulas $(a)$ and $(b)$, we have
\[
(c)\;\;\;\;\;\;\;\sum_{p=1}^{2i+1}x_p\binom{2j}{p}=j\left[\binom{2j}{2i}+\binom{2K-2j}{2i}\right]
\]
and
\[
\begin{split}
(d)\;\;\;&\;\;\;\sum_{p=1}^{2i+1}x_p\binom{2j}{2K-p-1}\\
&=\sum_{p=1}^{2i+1}x_p\binom{2j}{2j-2K+p+1}\\
&=\sum_{p=0}^{2i+1}x_{p+(2K-2j-1)}\binom{2j}{p}\\
&=j\binom{2j}{2i+2j-2K+1}+\frac{2K-2j-1}{2}\binom{2j+1}{2i+2j-2K+2}\\
&\;\;-j\binom{1}{2i+2j-2K+1}+\frac{2k-2j-1}{2}\binom{1}{2i+2j-2K+2}\\
&=j\binom{2j}{2K-2i-1}+\frac{2K-2j-1}{2}\binom{2j+1}{2K-2i-1}\\
&\;\;-j\binom{1}{2i+2j-2K+1}+\frac{2k-2j-1}{2}\binom{1}{2i+2j-2K+2}.
\end{split}
\]

By formulas $(c)$ and $(d)$,
\[\begin{split}
(e)\;\;\;\;\;\sum_{l=1}^{i+1}x_{2l-1}&\left[\binom{2j}{2K-2l}+\binom{2j}{2l-1}-\delta\binom{j}{K-l}\right]\\
                &+\sum_{l=1}^{i}x_{2l}\left[\binom{2j}{2l}+\binom{2j}{2K-2l-1}-\delta\binom{j}{l}\right]
\end{split}\]
\[
=\sum_{p=1}^{2i+1}x_p\binom{2j}{p}+\sum_{p=1}^{2i+1}x_p\binom{2j}{2K-p-1}-\sum_{l=1}^{i+1}x_{2l-1}\delta\binom{j}{K-l}-\sum_{l=1}^{i}x_{2l}\delta\binom{j}{l}
\]
\[
\begin{split}
=&j\left[\binom{2j}{2i}+\binom{2K-2j}{2i}\right]+j\binom{2j}{2K-2i-1}+\frac{2K-2j-1}{2}\binom{2j+1}{2K-2i-1}\\
 &-j\binom{1}{2i+2j-2K+1}+\frac{2K-2j-1}{2}\binom{1}{2i+2j-2K+2}\\
 &-\sum_{l=1}^{i+1}x_{2l-1}\delta\binom{j}{K-l}
 -\sum_{l=1}^{i}x_{2l}\delta\binom{j}{l}\\
=&j\left[\binom{2j}{2i}+\binom{2j}{2K-2i-1}-\delta\binom{j}{i}+\binom{2K-2j}{2i}+\binom{2K-2j}{2K-2i-1}-\delta\binom{K-j}{i}  \right].
\end{split}
\]
While the last identity follows from Lemma \ref{ls} below.
From the formula $(e)$, we have
\[
\begin{split}
\sum_{j=1}^{K-1}jc_j&\left[\binom{2j}{2i}+\binom{2j}{2K-2i-1}-\delta\binom{j}{i}\right]+\\
                    &\sum_{j=1}^{K-1}jc_j\left[\binom{2K-2j}{2i}+\binom{2K-2j}{2K-2i-1}-\delta\binom{K-j}{i}\right]
\end{split}
\]
\[\begin{split}
&=\sum_{l=1}^{i+1}x_{2l-1}\sum_{j=1}^{K-1}c_j\left[\binom{2j}{2K-2l}+\binom{2j}{2l-1}-\delta\binom{j}{K-l} \right]\\
&\;\;\;\;\;+\sum_{l=1}^ix_{2l}\sum_{j=1}^{K-1}c_j\left[\binom{2j}{2l}+\binom{2j}{2K-2l-1}-\delta\binom{j}{l}  \right]\\
&=0.
\end{split}\]
Thus the lemma is proved.    $\hfill\Box$\\

\begin{lem}\label{ls}
Let $\{{x_l}\}$ be as above. For $1\le i,j \le K-1$ we have
\begin{align*}
&-\sum_{l=1}^{i+1}x_{2l-1}\delta \binom{j}{K-l}-\sum_{l=1}^ix_{2l}\delta \binom{j}{l}\\
=&j\left [\binom{1}{2i+2j-2K+1}-\delta \binom{j}{i}+\binom{2K-2j}{2K-2i-1} -\delta \binom{K-j}{i} \right ]\\
& -\frac{2K-2j-1}{2}\left[\binom{1}{2i+2j-2K+2}+\binom{2j+1}{2K-2i-1}\right].
\end{align*}
\end{lem}

\noindent{\bf Proof:} 
The proof consists of checking 9 cases separately. For convenience, $LHS$ and $RHS$ are short for left hand side and right hand side respectively.\\

\textbf{Case I:} $i<K-1-i$, i.e. $K>2i+1$.
\begin{itemize}
\item [Subcase 1:] $1\le j \le i$.
\begin{align*}
LHS &=-x_{2j} \\
&=-j\binom{1}{2i-2j+1}+j\binom{2K-2j}{2i-2j+1} \\
 &=-j\binom{1}{2i-2j+1}+j\binom{1}{2i+2j-2K+1}+j\binom{2K-2j}{2i-2j+1}\\
 &=-j\delta\binom{j}{i}+j\binom{1}{2i+2j-2K+1}+j\binom{2K-2j}{2i-2j+1}\\ &=RHS.
 \end{align*}
 
 \item[Subcase 2:] $i<j <K-1-i.$ $LHS=0=RHS$.
 
 \item[Subcase 3:] $K-1-i\le j \le K-1$.
 \begin{align*}
 LHS &=-x_{2K-2j-1}\\
        &=-\frac{2K-2j-1}{2}\left[\binom{1}{2i+2j-2K+2} +\binom{2j+1}{2i+2j-2K+2} \right]\\
        & \quad (\text{denote the above expression by}-M)\\
        &=RHS.
 \end{align*}
\end{itemize} 

\textbf{Case II:} $i=K-i-1$, i.e. $K=2i+1$.
\begin{itemize}
\item [Subcase 1:] $1\le j < i$.
\begin{align*}
LHS &=-x_{2j}=-j\binom{1}{2i-2j+1}+j\binom{2K-2j}{2i-2j+1} \\
& =j\binom{2K-2j}{2i-2j+1}=RHS.
\end{align*}

\item [Subcase 2:] $j=i$.
\[
LHS=-x_{2i}-x_{2i+1}=(2i+1)(i-1)=RHS.
\]

\item[Subcase 3:] $i<j\le K-1$.
\begin{align*}
LHS=-x_{2K-2j-1} &=-M \\
& =-\frac{2K-2j-1}{2} \binom{2j+1}{2i+2j-2K+2}. 
\end{align*}
\begin{align*}
RHS&=-j\left[\delta \binom{K-j}{i}-\binom{1}{2i+2j-2K+1}\right]  \\
& \quad -\frac{2K-2j-1}{2} \binom{2j+1}{2i+2j-2K+2}\\
& =-\frac{2K-2j-1}{2} \binom{2j+1}{2i+2j-2K+2}=LHS.
\end{align*}
\end{itemize}

\textbf{Case III:} $K-i-1<i,$ i.e. $K<2i+1$.
\begin{itemize}
\item [Subcase 1:] $1\le j <K-i-1$. 
\begin{align*}
LHS=-x_{2j}=j\binom{2K+2i-2j-2}{2i-2j+1}=RHS.
\end{align*}

\item [Subcase 2:] $K-i-1 \le j \le i$.
\begin{align*}
LHS &=-x_{2j}-x_{2K-2j-1} \\
 &=-j\left [ \binom{1}{2i-2j+1}- \binom{2K-2j}{2i-2j+1}\right]-M \\
 &= -j \left[\delta\binom{j}{i}-\binom{2K-2j}{2i-2j+1}\right]-M\\
 &=RHS.
\end{align*}

\item [Subcase 3:] $i<j\le K-1$.
\begin{align*}
LHS &= -x_{2K-2j-1}=-M=-\frac{2K-2j-1}{2} \binom{2j+1}{2i+2j-2K+2} \\ 
&=RHS.
\end{align*}
\end{itemize} $\hfill\Box$\\

By Lemma \ref{odd}, \ref{ev} and \ref{sym}, now we can prove Proposition \ref{poly}.

\noindent{\bf Proof of Proposition \ref{poly}}:
For $C(X)=\sum\limits_{n=1}^{\frac{1}{2}(N-3)}c_{2n} X^{2n}$, let
\[
p(X)=C\;'(X)+X^{N-3}C\;'(\frac{1}{X})
\]
and $q(X)=X^{N-1}C(\frac{1}{X})-C(X)$.

It is easy to check that
\[
p(X)=X^{N-3}p(\frac{1}{X}),\;q(X)+X^{N-1}q(\frac{1}{X})=0,
\]
\[
X^{N-2}q\;'(\frac{1}{X})=(N-1)C(X)-Xp(X),\;q(0)=0.
\]
By Lemma \ref{sym}, if $C(X)$ satisfies $L_C(X)=\mathrm{constant}+\mathrm{odd\;polynomial\;in}\;X$, then $L_{Xp(X)}(X)=\mathrm{constant}+\mathrm{odd\;polynomial\;in}\;X$. By Lemma \ref{odd}, $p\in W_{N-1}^+$.

Since $$L_C(X)=\mathrm{constant}+\mathrm{odd\;polynomial\;in}\;X$$ and $$L_{Xp(X)}(X)=\mathrm{constant}+\mathrm{odd\;polynomial\;in}\;X,$$  we obtain
\[
L_{X^{N-2}q\;'(\frac{1}{X})}(X)=\mathrm{constant}+\mathrm{odd\;polynomial\;in}\;X.
\]
By Lemma \ref{ev}, $q\in W_{N+1}^-$.
$\hfill\Box$\\

Now from Proposition \ref{injec},  \ref{surj} and  \ref{poly}, Theorem \ref{exact} is proved. By counting dimensions of spaces in the exact sequence in Theorem \ref{exact}, Theorem \ref{isom} is proved. Thus the classification of all the relations among depth-graded motivic double zeta values $\zeta^{\mathfrak{m}}(r,N-r)$ with odd $r\geq 3$ and $N-r\geq 2$ is completed.
\begin{rem}\label{formal}
Theorem \ref{isom} is also true for the formal double zeta values. We omit the proof here.  For more  exact sequences about depth-graded motivic multiple zeta values, see \cite{li}.
\end{rem}

\section*{Acknowledgements}
The authors would like to express their sincere gratitude to the anonymous referee for  his/her detailed comments to improve this paper.

\end{document}